\documentclass[12pt]{article} 
\def\R{\hbox{{\rm I}\kern-0.2em{\rm R}\kern0.2em}}
\def\D{\hbox{{\rm I}\kern-0.2em{\rm D}\kern0.2em}}

\def\be{\begin{equation}}
\def\ee{\end{equation}}
 
\def\p{\partial}
\def\({\left(}
\def\){\right)}
\def\[{\left[}
\def\]{\right]}
\def\bc{\begin{center}}
\def\ec{\end{center}}

\begin{document}

{\large \bf Invariant Linearization Criteria for Systems of Cubically
Semi-Linear Second-Order Ordinary Differential Equations}

F M Mahomed$^1$ and Asghar Qadir$^{2}$

$^1$Centre for Differential Equations, Continuum Mechanics and Applications\\
School of Computational and Applied Mathematics\\
University of the Witwatersrand\\
Wits 2050, South Africa\\
E-mail: Fazal.Mahomed@wits.ac.za

$^2$Centre for Advanced Mathematics and Physics\\
National University of Sciences and Technology\\
Campus of the College of Electrical and Mechanical Engineering\\
Peshawar Road, Rawalpindi, Pakistan
and\\
Department of Mathematical Sciences\\
King Fahd University of Petroleum and Minerals\\
Dhahran 31261, Saudi Arabia\\
E-mails: aqadirs@comsats.net.pk

{\bf Abstract}. Invariant linearization criteria of square systems of
second-order quadratically semi-linear ordinary differential
equations (ODEs) that can be represented as geodesic equations are
extended to square systems of ODEs cubically nonlinear in the first derivatives. 
It is shown that there are two branches for the linearization problem via point transformations
for an arbitrary system of second-order
ODEs. One is when the system is at most cubic in the first derivatives. 
We solve this branch of the linearization problem by point transformations 
in the case of a square sytem of two second-order ODEs.
Necessary and sufficient conditions
for linearization by means of point transformations
are given in terms of coefficient functions of the system 
of two second-order ODEs cubically nonlinear in the first derivatives.
A consequence of our geometric approach of projection is a
re-derivation of Lie's conditions for a single second-order ODE and
sheds light on more recent results on them. In particular, we show here
how one can construct point transformations for reduction to the
simplest linear equation by going to the higher space and just utilising
the coefficients of the original ODE. We also obtain invariant criteria
for the reduction of a linear square system to the simplest system.
Moreover, these results contain the quadratic case as a special case. Examples are
given to illustrate our results.

\section{Introduction}

A linearization problem involves the study of families of equations
that are reducible via admissible transformations, which can be
point, contact or more general, to linear equations. Lie \cite{lie}
presented  linearizability criteria, obtaining both algebraic and
practical criteria, for a single second-order ODE to be point
transformable to a linear equation via invertible changes of both
the independent and dependent variables.

Lie \cite{lie} proved that necessary and sufficient conditions for a
second-order ODE, $y''=E(x,y,y')$, to be linearizable by means of
invertible point transformations are that the ODE be at most cubic
in the first derivative, viz.
\begin{equation}
y''+E_3(x,y)y'^3+E_2(x,y)y'^2+E_1(x,y)y'+E_0(x,y)=0 \label{(1)}
\end{equation}
and the coefficients $E_0$ to $E_3$ satisfy the over-determined
integrable system
\begin{eqnarray}
b_x&=&-\frac13 E_{1y}+\frac23 E_{2x}+be-E_0E_3,\nonumber\\
b_y&=&E_{3x}-b^2+bE_2-E_1E_3+eE_3,\nonumber\\
e_x&=&E_{0y}+e^2-eE_1-bE_0+E_0E_2,\nonumber\\
e_y&=&\frac23E_{1y}-\frac13E_{2x}-be+E_0E_3,\label{(2)}
\end{eqnarray}
where $b$ and $e$ are auxiliary variables and the suffices $x$ and
$y$ here and hereafter refer to partial derivatives. Since the
classic work of Lie there has been continuing interest in this
topic. We, inter alia, re-derive the Lie conditions (\ref{(2)})
geometrically, by projections.

Tress\'e \cite{tre} also studied the linearization problem for
scalar second-order ODEs. He deduced two relative invariants of the
equivalence group of point transformations, the vanishing of both of
which gives necessary and sufficient conditions for linearization of
equation (\ref{(1)}). These conditions are equivalent to the Lie
conditions (2) (see Mahomed and Leach \cite{mah1}) and can be given
as the compatibility of (2) as
\begin{eqnarray}
3(E_1E_3)_x-E_{1yy}+2E_{2xy}-3(E_0E_3)_y+E_2E_{1y}-2E_2E_{2x}-3E_{3xx}-3E_3E_{0y}=0,\nonumber\\
3(E_0E_3)_x+2E_{1xy}-3E_{0yy}-E_{2xx}-E_1E_{2x}+2E_1E_{1y}-3(E_0E_2)_y+3E_0E_{3x}=0.\label{(3)}
\end{eqnarray}
Note that under the interchange of $E_3$ by $-E_0$, $E_2$ by $-E_1$
and $x$ by $y$, these conditions imply each other. Equations
(\ref{(3)}) provide practical criteria for linearization of equation
(\ref{(1)}) by point transformations. These conditions were also
derived by the Cartan equivalence method (see Grissom et al
\cite{gri}) as well as recently using a geometric argument in
Ibragimov and Magri \cite{mag}. The reader is also referred to the review of various approaches in Mahomed
\cite{mah}. Linearization via point and other
than point transformations is of great interest and has been
investigated in several works (see, e.g.
\cite{che,mah2,mah3,gre,waf1,neu,mel1,mel2}).

The algebraic criteria of linearization of systems of second-order
ODEs by means of point transformations have been considered in Wafo
and Mahomed \cite{waf1}. Practical criteria for quadratic
semi-linear systems of second-order ODEs have been looked at
recently as well (see Mahomed and Qadir \cite{mah4}). In this paper
our intention is to extend these results to cubically semi-linear
square systems of second-order ODEs using geometric methods
developed earlier (see Feroze et al \cite{fer}). As a by-product of
our approach we re-derive the Lie conditions (\ref{(2)}). Moreover,
we present practical criteria in terms of coefficients for cubically
semi-linear systems of second-order ODEs to be linearizable by point
transformations. As a consequence we provide practical criteria for
the class of linear second-order system of two ODEs to be reducible
to the simplest system. Notwithstanding, our results subsume the
linearization criteria for the quadratic case.

The outline of this paper is as follows. In the next section we
present mathematical preliminaries. In section 3 we give an
alternative method for obtaining the Lie conditions (\ref{(2)}) as
well as an alternative method for the construction of linearizing
transformations for scalar second-order ODEs. Then in section 4 we
derive practical criteria for linearization for a system of two
second-order cubically semi-linear ODEs. Herein we state the
relevant result for linear systems. Our theorem also contains the
quadratically semi-linear equations as a corollary. In the next section we
provide examples that amply illustrate our results. Finally, in section 6 we 
present a brief summary and conclusion.

\section{Preliminaries}

We first present some preliminaries. The system of geodesic
equations is
\begin{equation}
{\ddot x}^i+\Gamma^i_{jk}{\dot x}^j{\dot x}^k=0, \quad
i,j,k=1,\ldots,n,\label{(4)}
\end{equation}
where the dot refers to total differentiation with respect to the
parameter $s$ and $\Gamma^i_{jk}$ are the Christoffel symbols, which
depend on $x^i$ and are given in terms of the metric tensor as
\begin{equation}
\Gamma^i_{jk}=\frac12g^{im}(g_{jm,k}+g_{km,j}-g_{jk,m}).\label{(5)}
\end{equation}
The Christoffel symbols are symmetric in the lower pair of indices
and have $n^2(n+1)/2$ coefficients. The Riemann curvature tensor is
\begin{equation}
R^i_{\,jkl}=
\Gamma^i_{jl,k}-\Gamma^i_{jk,l}+\Gamma^i_{mk}\Gamma^m_{jl}-\Gamma^i_{ml}\Gamma^m_{jk},\label{(6)}
\end{equation}
which is skew-symmetric in the lower last two indices and satisfies
\begin{equation}
R^i_{\,jkl}+R^i_{\,klj}+R^i_{\,ljk}=0.\label{(7)}
\end{equation}
A necessary and sufficient condition for a system of $n$
second-order quadratically semi-linear ODEs for $n$ dependent
variables of the form (\ref{(4)}) to be linearizable by point
transformation and admit $sl(n+2,\R)$ symmetry algebra is that the
Riemann tensor vanishes (\cite{mah3,ami}), i.e.
\begin{equation}
R^i_{\,jkl}=0.\label{(8)}
\end{equation}
Practical criteria and the construction of point transformations are
given in \cite{mah4}. In particular, for a system of two geodesic
equations (\ref{(4)}), one has the linearization conditions
(admittance of $sl(4,\R)$ symmetry algebra) on the coefficients
given by
\begin{eqnarray}
a_y-b_x + be -cd &=& 0,\nonumber\\
b_y -c_x+ (ac - b^2) + (bf - ce) &=& 0,\nonumber\\
d_y-e_x - (ae - bd) - (df - e^2) &=& 0,\nonumber\\
(b + f)_x &=& (a + e)_y,\label{(9)}
\end{eqnarray}
where the Christoffel symbols are
\begin{equation}
\Gamma^1_{11} = -a, \Gamma^1_{12} = -b, \Gamma^1_{22} = -c,
\Gamma^2_{11} = -d, \Gamma^2_{12} = -e, \Gamma^2_{22} =
-f.\label{(10)}
\end{equation}
Now equation (\ref{(5)}) together with (\ref{(10)}) on setting
$g_{11}=p$, $g_{12}=q=g_{21}$ and $g_{22}=r$ yield
\begin{eqnarray}
p_x &=& -2(ap + dq),\nonumber\\
q_x &=& -bp - (a + e)q - dr,\nonumber\\
r_x &=& -2(bq + er),\nonumber\\
p_y &=& -2(bp + eq),\nonumber\\
q_y &=& -cp - (b + f)q - er,\nonumber\\
r_y &=& -2(cq + fr).\label{(11)}
\end{eqnarray}

The construction of the linearization point transformations are
found as follows (see \cite{mah4}). One invokes
\begin{equation}
g_{ab}({\bf x})={\p u^i\over\p x^a}{\p u^j\over\p x^b}g_{ij}({\bf
u}),\label{(12)}
\end{equation}
where ${\bf x}=(x^1,\ldots,x^n)$, ${\bf u}=(u^1,\ldots,u^n)$ with
the requirement that $g_{ij}({\bf u})$ be the identity matrix. For
the case of two variables, we need to solve the equations
\begin{equation}
u_x^2 + v_x^2 = p,\; u_x u_y + v_x v_y = q,\; u_y^2 + v_y^2 =
r,\label{(13)}
\end{equation}
for which we have set $(x^1,x^2)=(x,y)$, $(u^1,u^2)=(u,v)$,
$g_{11}=p$, $g_{12}=q=g_{21}$ and $g_{22}=r$ in (\ref{(12)}).

Following Aminova and Aminov \cite{ami}, we project the system down
by one dimension and write the geodesic equations (\ref{(4)}) as
\begin{equation}
{x^a}''+A_{bc}{x^a}'{x^b}'{x^c}'+B^a_{bc}{
x^b}'{x^c}'+C^a_b{x^b}'+D^a=0,\quad a=2,\ldots,n,\label{(14)}
\end{equation}
where the prime now denotes differentiation with respect to the
parameter $x^1$ (in \cite{ami} $x^n$ is used as the parameter) and
the coefficients in terms of the $\Gamma^a_{bc}$s are
\begin{equation}
A_{bc}=-\Gamma^1_{bc},\,
B^a_{bc}=\Gamma^a_{bc}-2\delta^a_{(c}\Gamma^1_{b)1},\,
C^a_b=2\Gamma^a_{1b}-\delta^a_b\Gamma^1_{11},\, D^a=\Gamma^a_{11},\;
a,b,c=2,\ldots,n,\label{(15)}
\end{equation}
where we have used the notation $T_{(a,b)}=(T_{ab}+T_{ba})/2$.
It is straightforward to deduce (\ref{(14)}) and (\ref{(15)}).
Indeed, insert
$$
{\dot x}^a={d x^a\over d x^1}{\dot x}^1,\quad a=2,\ldots,n
$$
and its derivatives
$$
{\ddot x}^a={d^2 x^a\over d {x^1}^2}\dot{x^1}^2+ {d x^a\over
dx^1}{\ddot x}^1,\quad a=2,\ldots,n
$$
 into system (\ref{(4)}). These, after cancelation of $\dot{x^1}^2$, directly yield
(\ref{(14)}) and (\ref{(15)}). Note that in projecting down the
Christoffel symbols there is degeneracy which results from the
reduction of the range of the indices, so that $\Gamma_{11}^1$ and
$\Gamma_{b1}^1$ appear in the same combinations in $C^a_b$ and
$B^a_{bc}$, respectively. Consequently the set of coefficients ${\bf
A}$, ${\bf B}$, ${\bf C}$, ${\bf D}$ have $n$ less elements than the
coefficients $\Gamma^i_{jk}$.

\section{Re-derivation of the Lie conditions}

We invoke equations (\ref{(14)}) and (\ref{(15)}) for $n=2$. We also
use (\ref{(10)}) in identifying the $\Gamma^i_{jk}$s with the
coefficients $a$ to $f$ of the system of two geodesic equations
which projects to (\ref{(14)}). Thus we have (setting $(x^1,x^2)=(x,y)$)
\begin{equation}
y''+E_3(x,y)y'^3+E_2(x,y)y'^2+E_1(x,y)y'+E_0(x,y)=0,\label{(16)}
\end{equation}
where
\begin{eqnarray}
E_3&=&A_{22}=-\Gamma_{22}^1=c,\nonumber\\
E_2&=&B^2_{22}=\Gamma_{22}^2-2\Gamma_{12}^1=-f+2b,\nonumber\\
E_1&=&C^2_2=2\Gamma_{12}^1-\Gamma_{11}^1=-2e+a,\nonumber\\
E_0&=&D^2=\Gamma_{11}^2=-d.\label{(17)}
\end{eqnarray}
To re-derive the Lie conditions (\ref{(2)}), we use the system of
two geodesic equations (\ref{(4)}) from which equation (\ref{(16)})
arises projectively. Hence we utilize the conditions (\ref{(9)})
which are conditions for a flat space. This requires that the
coefficients $a$ to $f$ be in terms of the $E_i$s. From (\ref{(17)})
we have
\begin{eqnarray}
a&=&E_1+2e,\nonumber\\
c&=&E_3,\nonumber\\
d&=&-E_0,\nonumber\\
f&=&2b-E_2, \label{(18)}
\end{eqnarray}
where we have chosen $b$ and $e$ as yet arbitrary. These are
constrained by the relations (\ref{(9)}). We substitute (\ref{(18)})
into (\ref{(9)}). Equations (\ref{(9)}) then yield
\begin{eqnarray}
E_{1y}+2e_y-b_x+be+E_3E_0=0,\nonumber\\
b_y-E_{3x}+E_1E_3+eE_3+b^2-bE_2=0,\nonumber\\
E_{0y}+e_x+eE_1+e^2-bE_0+E_0E_2=0,\nonumber\\
3b_x-3e_y-E_{2x}-E_{1y}=0.\label{(19)}
\end{eqnarray}
The first and last equations of (\ref{(19)}) are easily seen to be
equivalent to
\begin{eqnarray}
b_x&=&-\frac13 E_{1y}+\frac23E_{2x}-be-E_0E_3,\nonumber\\
e_y&=&\frac13E_{2x}-\frac23E_{1y}-be-E_0E_3.\label{(20)}
\end{eqnarray}
The second and third equations of (\ref{(19)}) as well as equations
(\ref{(20)}), on replacing $e$ by $-e$, are precisely the Lie
conditions (\ref{(2)}). Hence, we have provided an alternative
derivation of the Lie conditions (\ref{(2)}) by viewing the
projection (\ref{(16)}) in one higher space and looking at the flat
space requirement there. If we had projected the system of two
geodesic equations to a single ODE of the form (\ref{(16)}) by using
$x^2$ instead of $x^1$, then by interchanging $E_3$ by $-E_0$, $E_2$
by $-E_1$ and $x^1=x$ by $x^2=y$, the coefficients (\ref{(17)})
imply the coefficients of the projected equation with independent
variable $x^2$. We state the following theorem.

{\bf Theorem 1}. A necessary and sufficient condition that the
scalar second-order ODE (\ref{(16)}) has $sl(3,\R)$ symmetry algebra
is that there is a
corresponding system of two geodesic equations of the form
(\ref{(4)}) from which it is projected that admits the $sl(4,\R)$ symmetry
algebra.

Furthermore, one can construct linearizing point transformations
for (\ref{(16)}) that satisfy (\ref{(3)}) by resorting to the
corresponding system of two geodesic equations from which
(\ref{(16)}) arises by projection. This is done by using the
relations (\ref{(13)}). This approach also results in the
determination of at least one metric as a bonus. Notwithstanding, this method uses the coefficients of 
the equation which is linearizable and a transformation is then constructed via the relations (\ref{(13)}).
We consider two examples to illustrate this.

{\bf 1.} On using (\ref{(18)}), the simple nonlinear equation
\begin{equation}
y''+y'^3-y'=0\label{(21)}
\end{equation}
has corresponding $a$ to $f$ values, 
$$
a=-1+2e, c=1, d=0, f=2b.
$$
These together with the choices $b=0$ and $e=1$ satisfy the system
(\ref{(9)}). With these values of $a$ to $f$ we obtain from
(\ref{(11)}) particular solutions for $p$, $q$ and $r$ given by
$$
p=r=\exp(2y-2x),\quad q=-\exp(2y-2x).
$$
Invoking (\ref{(13)}), a linearizing point transformation to the
simplest second-order ODE is
$$
u=\frac1{\sqrt{2}}\exp(-x+y),\quad v=\frac1{\sqrt{2}}\exp(-x-y),
$$
where $u$ is the new independent variable.

{\bf 2.} The familiar nonlinear ODE (see, e.g. \cite{lea1})
$$
y''+3yy'+y^3=0
$$
has, upon using (\ref{(18)}),
$$
a=3y+2e, c=0, d=-y^3, f=2b.
$$
These and the choices $b=1/y$ and $e=-y$ satisfy (\ref{(9)}). A
particular solution of (\ref{(11)}) is then
$$
p=1+x^2-2xy^{-1}+y^{-2},\; q=(1+x^2)y^{-2}-xy^{-3},\;
r=y^{-4}(1+x^2).
$$
A  point transformation that linearizes the original ODE to the
simplest second-order equation, after solving (\ref{(13)}),  then is
$$
u=x-y^{-1},\quad v=\frac12 x^2-{x\over y},
$$
where $u$ is taken as the new independent variable.
This transformation was previously obtained in \cite{lea1} by
mapping generators to canonical forms. As such we have presented another
way of finding such transformations.

\section{Linearization conditions for square systems}

Driven by the success in obtaining the Lie conditions (\ref{(2)}) by
projection and then going back to the geodesic equations, we pursue
similar conditions and practical criteria for linearization for a
system of two second-order ODEs in a similar manner. Consequently, we study
(\ref{(14)}) for linearization via point transformations by
resorting to a system of three geodesic equations (\ref{(4)}).
Before we do so, we need to first understand what is meant by
linearization for systems of ODEs. A system of two second-order
linear ODEs can possess 5, 6, 7, 8 or 15 point symmetries (see
\cite{gor,waf2}). The maximal symmetry
algebra and hence $sl(4,\R)$ is reached for the simplest system. Here we consider
practical linearization criteria in terms of the coefficients for a
system of two cubically semi-linear second-order ODEs of the form
(\ref{(14)}) having $sl(4,\R)$ symmetry algebra. The quadratically
semi-linear case was treated in \cite{mah4}. Also algebraic criteria
for systems of second-order ODEs have been found in \cite{waf1}.

We once again invoke equations (\ref{(14)}) and (\ref{(15)}) but now
for $n=3$. We therefore have
\begin{eqnarray}
{x^2}''+A_{22}({x^2}')^3+2A_{23}({x^2}')^2{x^3}'+A_{33}{x^2}'({x^3}')^2+B^2_{22}({
x^2}')^2+2B^2_{23}{x^2}'{x^3}'\nonumber\\
+B_{33}^2({x^3}')^2+C^2_2{x^2}'+C^2_3{x^3}'+D^2=0,\nonumber\\
{x^3}''+A_{22}({x^2}')^2{x^3}'+2A_{23}{x^2}'({x^3}')^2+A_{33}({x^3}')^3+B^3_{22}({
x^2}')^2+2B^3_{23}{x^2}'{x^3}'\nonumber\\
+B_{33}^3({x^3}')^2+C^3_2{x^2}'+C^3_3{x^3}'+D^3=0,\label{(22)}
\end{eqnarray}
with coefficients
\begin{equation}
A_{bc}=-\Gamma^1_{bc},\,
B^a_{bc}=\Gamma^a_{bc}-2\delta^a_{(c}\Gamma^1_{b)1},\,
C^a_b=2\Gamma^a_{1b}-\delta^a_b\Gamma^1_{11},\, D^a=\Gamma^a_{11},\;
a,b,c=2,3.\label{(23)}
\end{equation}

Here three $\Gamma_{bc}^a$ coefficients are lost. We select
$\Gamma_{12}^1$, $\Gamma_{12}^2$ and $\Gamma_{33}^3$ as arbitrary.
We solve for the 15 $\Gamma^a_{bc}$s of (\ref{(23)}) in terms of the
15 coefficients  $A_{bc}$, $B^a_{bc}$, $C^a_b$, $D^a$ as well as
$\Gamma_{12}^1$, $\Gamma_{12}^2$ and $\Gamma_{33}^3$. We only write
down the $\Gamma^a_{bc}$s in which the arbitrary elements appear.
They are
\begin{eqnarray}
\Gamma_{11}^1&=&2\Gamma_{12}^2-C^2_2,\nonumber\\
\Gamma_{13}^1&=&\frac12(\Gamma_{33}^3-B^3_{33}),\nonumber\\
\Gamma_{22}^2&=&2\Gamma_{12}^1+B^2_{22},\nonumber\\
\Gamma_{23}^2&=&\frac12(\Gamma_{33}^3+2B^2_{23}-B^3_{33}),\nonumber\\
\Gamma_{13}^3&=&\Gamma_{12}^2+\frac12C^3_3-\frac12C^2_2,\nonumber\\
\Gamma_{23}^3&=&\Gamma_{12}^1+B^3_{23},\label{(24)}
\end{eqnarray}
The others can be read-off from equations (\ref{(23)}).

The flat space requirement for the corresponding system of three
geodesic equations (\ref{(4)}) are now imposed by means of the
vanishing of the Riemann tensor, viz. (\ref{(8)}). They are (let
$(x^1,x^2,x^3)=(x,y,z)$)
\begin{eqnarray}
(\Gamma_{j2}^i)_x-(\Gamma_{j1}^i)_y+\Gamma_{m1}^i\Gamma_{j2}^m-\Gamma_{m2}^i\Gamma_{j1}^m=0,
\nonumber\\
(\Gamma_{j3}^i)_x-(\Gamma_{j1}^i)_z+\Gamma_{m1}^i\Gamma_{j3}^m-\Gamma_{m3}^i\Gamma_{j1}^m=0,
\nonumber\\
(\Gamma_{j3}^i)_y-(\Gamma_{j2}^i)_z+\Gamma_{m2}^i\Gamma_{j3}^m-\Gamma_{m3}^i\Gamma_{j2}^m=0,
\label{(25)}
\end{eqnarray}
which provide 27 conditions. Only 24 of them are linearly
independent due to the identity (\ref{(7)}). The reduction of these
equations to explicit form is given in the Appendix.

These are 24 conditions (\ref{(A1)}) to (\ref{(A3)}) given in the
Appendix that arise from the vanishing of the Riemann tensor as
given in (\ref{(25)}). They are the Lie-type integrability
conditions for the $\Gamma^i_{jk}$. We find that there are 7
equations in (\ref{(A1)}) to (\ref{(A3)}) which are independent of
the $\Gamma^i_{jk}$. The other 17 contain first-order partial
derivatives of the $\Gamma^i_{jk}$. Of these, $\Gamma_{12,y}^2$ and
$\Gamma_{12,z}^2$ appear once each, $\Gamma_{33,x}^3$ occurs three times
and the rest twice each. Therefore, apart from the 7 conditions
which are independent of the $\Gamma^i_{jk}$ and given solely in
terms of the coefficients of the system, there arise a further 8
conditions on the coefficients upon equating the respective
$\Gamma^i_{jk}$. Hence, we end up with 15 conditions or constraint
equations on the coefficients. Now the $\Gamma^i_{jk}$ which appear
once each do not result in linearly independent equations as can
easily be checked by equating them with the corresponding
$\Gamma^i_{jk}$ that were discarded. The resultant two equations
that occur in this manner are linearly dependent. Thus the
$\Gamma_{12,y}^2$ and $\Gamma_{12,z}^2$ are spurious. It is thus
opportune to state the following theorem.

{\bf Theorem 2}. A necessary and sufficient condition for the system
of two cubically semi-linear ODEs
\begin{eqnarray}
&&y''+A_{22}y'^3+2A_{23}y'^2z'+A_{33}y'z'^2+B_{22}^2y'^2+2B_{23}^2y'z'+B_{33}^2z'^2+C^2_2y'+C^2_3z'+D^2=0,\nonumber\\
&&z''+A_{22}y'^2z'+2A_{23}y'z'^2+A_{33}z'^3+B_{22}^3y'^2+2B_{23}^3y'z'+B_{33}^3z'^2+C^3_2y'+C^3_3z'+D^3=0,\nonumber\\
\label{(50)}
\end{eqnarray}
(where the prime denotes differentiation with respect to the independent variable $x$
and the coefficients are in general functions of $x,y,z$) to be linearizable via point transformations to the
simplest system of two second-order ODEs is that its coefficients
satisfy the following fifteen conditions on the coefficients
functions of (\ref{(50)}), viz.
\begin{eqnarray}
&&\frac12C^3_{2x}-D^3_y+\frac14C^3_3C^3_2+\frac14C^2_2C^3_2-D^2B^3_{22}-D^3B^3_{23}=0,\nonumber\\
&&B^3_{22x}-\frac12C^3_{2y}-A_{22}D^3+\frac12C^3_2B^2_{22}+\frac12C^3_3B^3_{22}-\frac12C^2_2B^3_{22}-\frac12B^3_{23}C^3_2=0,
\nonumber\\
&&B^3_{23x}-\frac13B^2_{22x}+\frac16C^2_{2y}-\frac43D^3A_{23}-\frac23B^3_{22}C^2_3+\frac2{3}B^2_{23}C^3_2-\frac12C^3_{3y}=0,\nonumber\\
&&\frac12C^2_{3x}-D^2_z+\frac14C^2_3C^3_3+\frac14C^2_3C^2_2-B^2_{23}D^2-B^2_{33}D^3=0,\nonumber\\
&&B^2_{33x}-\frac12C^2_{3z}-D^2A_{33}+\frac12C^2_3B^3_{33}-\frac12B^2_{23}C^2_3-\frac12B^2_{33}C^3_3+\frac12B^2_{33}C^2_2=0,\nonumber\\
&&-A_{23y}+A_{22z}-A_{22}B^2_{23}-A_{23}B^3_{23}+A_{23}B^2_{22}+A_{33}B^3_{22}=0,\nonumber\\
&&-A_{33y}+A_{23z}-A_{22}B^2_{33}-A_{23}B^3_{33}+A_{23}B^2_{23}+A_{33}B^3_{23}=0,\nonumber\\
&&-A_{23x}+\frac5{6}A_{23}C^2_2+\frac13A_{33}C^3_2-\frac13B^3_{23z}+B^2_{33}B^3_{22}+\frac1{6}C^3_3A_{23}-B^2_{23}B^3_{23}\nonumber\\
&&-\frac23B^2_{23y}+\frac13B^3_{33y}+\frac23B^2_{22z}-\frac13C^2_3A_{22}=0,\nonumber\\
&&-A_{33x}+\frac12C^2_2A_{33}+\frac12A_{33}C^3_3-B^2_{33y}+B^2_{23z}-B^2_{22}B^2_{33}+B^2_{23}B^2_{23}-B^2_{23}B^3_{33}+
B^2_{33}B^3_{23}=0,\nonumber\\
&&-\frac23B^2_{22x}+\frac13C^2_{2y}-\frac12C^3_2B^3_{33}+D^2A_{22}-\frac23D^3A_{23}-\frac13C^2_3B^3_{22}+\frac5{6}B^2_{23}C^3_2\nonumber\\
&&+B^3_{23x}-\frac12C^3_{2z}+\frac12C^3_3B^3_{23}-\frac12C^2_2B^3_{23}=0,\nonumber\\
&&-A_{22x}+\frac12C^2_2A_{22}-B^3_{22}B^3_{33}+B^3_{23y}-B^3_{22z}+B^3_{22}B^2_{23}+B^3_{23}B^3_{23}+\frac12C^3_3A_{22}-
B^3_{23}B^2_{22}=0,\nonumber\\
&&D^2_y+B^2_{22}D^2+D^3B^2_{23}-D^3B^3_{33}+\frac12C^3_{3x}-\frac12C^2_{2x}-D^3_z+\frac14C^3_3C^3_3-\frac14C^2_2C^2_2-B^3_{23}D^2=0,
\nonumber\\
&&-2A_{23x}+\frac43B^3_{33y}+\frac13A_{23}C^2_2+\frac53A_{23}C^3_3+\frac23C^2_3A_{22}-\frac43B^3_{23z}-\frac23C^3_2A_{33}+2B^3_{22}B^2_{33}\nonumber\\
&&-2B^3_{23}B^2_{23}-\frac23B^2_{23y}+\frac23B^2_{22z}=0,\nonumber\\
&&B^2_{23x}+\frac12C^2_{3y}-2D^2A_{23}+\frac12C^2_3B^3_{23}+\frac12C^2_3B^2_{22}+\frac12C^3_3B^2_{23}-\frac12B^2_{23}C^2_2-B^2_{33}C^3_2\nonumber\\
&&-C^2_{2z}-D^3A_{33}=0,\nonumber\\
&&-B^2_{23x}+B^3_{33x}+C^2_{3y}-C^2_3B^3_{23}+C^2_3B^2_{22}+B^2_{23}C^3_3-B^2_{23}C^2_2-\frac12C^3_{3z}\nonumber\\
&&-\frac12C^2_{2z}-2D^3A_{33}=0.\label{(51)}
\end{eqnarray}

{\bf Proof}. The proof follows from the preceding discussions. For
if the system of two equations (\ref{(50)}) are linearizable by
point transformation to the simplest system, then its coefficients
can be written in terms of $\Gamma_{jk}^i$ as in equations
(\ref{(23)}) which in turn gives rise to the Lie-type integrability
conditions on the $\Gamma_{jk}^i$ and hence (\ref{(51)}).
Conversely, if the coefficients of the system of equations
(\ref{(50)}) satisfy the fifteen constraint conditions on the
coefficients given by the relations (\ref{(51)}) which is a
consequence of the Lie-type conditions (\ref{(A1)}) to (\ref{(A3)}),
then the coefficients of the system (\ref{(50)}) can be written in
terms of the $\Gamma_{jk}^i$ and the corresponding geodesic
equations in three-space is linearizable as well as the projected
equations (\ref{(50)}).

{\bf Corollary 1}. The system of two quadratically semi-linear ODEs
\begin{eqnarray}
&&y''+B_{22}^2y'^2+2B_{23}^2y'z'+B_{33}^2z'^2=0,\nonumber\\
&&z''+B_{22}^3y'^2+2B_{23}^3y'z'+B_{33}^3z'^2=0,
\label{(52)}
\end{eqnarray}
where the $B_{bc}^a$s are functions of $y$ and $z$ and the dot denotes total
derivative with respect to $x$, is linearizable by point
transformations to the simplest system of two equations if and only
if the $B^a_{bc}$s satisfy the four conditions on the coefficients
given by
\begin{eqnarray}
-B^3_{22}B^3_{33}+B^3_{23y}-B^3_{22z}+B^3_{22}B^2_{23}+B^3_{23}B^3_{23}-
B^3_{23}B^2_{22}&=&0,\nonumber\\
\frac43B^3_{33y}-\frac43B^3_{23z}+2B^3_{22}B^2_{33}
-2B^3_{23}B^2_{23}-\frac23B^2_{23y}+\frac23B^2_{22z}&=&0\nonumber\\
-\frac13B^3_{23z}+B^2_{33}B^3_{22}-B^2_{23}B^3_{23}
-\frac23B^2_{23y}+\frac13B^3_{33y}+\frac23B^2_{22z}&=&0,\nonumber\\
-B^2_{33y}+B^2_{23z}-B^2_{22}B^2_{33}+B^2_{23}B^2_{23}-B^2_{23}B^3_{33}+
B^2_{33}B^3_{23}&=&0.\label{(53)}
\end{eqnarray}

{\bf Remark}. If one sets $B^2_{22}=-a$, $B^2_{23}=-b$,
$B^2_{33}=-c$, $B^3_{22}=-d$, $B^3_{23}=-e$ and $B^3_{33}=-f$, one
gets precisely the conditions (\ref{(9)}). Hence Theorem 2 naturally
contains the linearizability criteria for the quadratic case.

{\bf Corollary 2}. The system of two linear (in the first
derivatives) ODEs
\begin{eqnarray}
&&y''+C^2_2y'+C^2_3z'+D^2=0,\nonumber\\
&&z''+C^3_2y'+C^3_3z'+D^3=0,\label{(54)}
\end{eqnarray}
where the prime refers to differentiation with respect to $x$ and the $C^a_b$s are
independent of $y$ and $z$, is linearizable by point transformations
to the simplest system of two equations if and only if the $C^a_b$s
and $D^a$s satisfy the three conditions on the coefficients, viz.
\begin{eqnarray}
\frac12C^3_{2x}+\frac14C^3_3C^3_2+\frac14C^2_2C^3_2&=&D^3_y,\nonumber\\
\frac12C^2_{3x}+\frac14C^2_3C^3_3+\frac14C^2_3C^2_2&=&D^2_z,\nonumber\\
\frac12C^3_{3x}-\frac12C^2_{2x}+\frac14C^3_3C^3_3-\frac14C^2_2C^2_2&=&D^3_z-D^2_y.
\label{(55)}
\end{eqnarray}

We have provided practical criteria, necessary and sufficient conditions, for equations of the form (\ref{(50)})
to be linearizable via point transformations to the simplest system. The question naturally arises if there are more general 
equations than (\ref{(50)}) that can be linearizable to the simplest system. Indeed, there are more general systems 
of two second-order ODEs which can be linearized.

The most general system of $n-1$ second-order ODEs linearizable is given by
\begin{equation}
J^i_j{x^j}''+G^i_{kj}{x^k}'{x^j}''+\Delta^i_{jkl}{x^j}'{x^k}'{x^l}'+
\Lambda^i_{jk}{x^j}'{x^k}'+\Omega^i_j{x^j}'+E^i=0, \; i=2,\ldots,n, \label{(r1)}
\end{equation}
where the prime refers to total differentiation with respect to $x^1$, the coefficient functions are dependent upon 
$x^1,\ldots,x^n$, and are given by
\begin{eqnarray}
&&J^i_j=X^1_{,1}X^i_{,j}-X^1_{,j}X^i_{,1},\nonumber\\
&&G^i_{kj}=X^1_{,k}X^i_{,j}-X^1_{,j}X^i_{,k},\nonumber\\
&&\Delta^i_{jkl}=X^1_{,l}X^i_{,jk}-X^1_{,jk}X^i_{,l},\nonumber\\
&&\Lambda^i_{jl}=2X^1_{,l}X^i_{,1j}-2X^1_{,1j}X^i_{,l}+X^1_{,1}X^i_{,jl}-X^i_{,1}X^1_{,jl},\nonumber\\
&&\Omega^i_j=2X^1_{,1}X^i_{,1j}-2X^1_{,1j}X^i_{,1}+X^1_{,j}X^i_{,11}-X^1_{,11}X^i_{,j},\nonumber\\
&&E^i=X^1_{,1}X^i_{,11}-X^1_{,11}X^i_{,1},\qquad i,j,k,l=2,\ldots,n
\label{(r2)}
\end{eqnarray}
in which 
\begin{equation}
X^1=X^1(x^1,\ldots, x^n), \quad X^i=X^i(x^1\ldots, x^n),\quad i=2,\ldots,n \label{(r3)}
\end{equation}
are invertible transformations. It is certainly not difficult to obtain (\ref{(r1)}). This is done by the substitution of
(\ref{(r3)}) into the free particle system 
\begin{equation}
{X^i}''=0,\quad i=2,\ldots,n; \quad '={d\over dX^1}.
\label{(r4)}
\end{equation}
This after routine calculations yields (\ref{(r1)}) with the coefficients satisfying (\ref{(r2)}).
Equation (\ref{(r1)}) is the most general system of $n-1$ equations point transformable to the 
simplest system (\ref{(r4)}). Equation (\ref{(r1)}) has $n(n-1)(n^2+6n-1)/6$ coefficients.

Equation (\ref{(r1)}) can be written in normal form in terms of at most cubic first order derivatives as
\begin{equation}
{x^i}''+A^i_{jkl}{x^j}'{x^k}'{x^l}'+B^i_{jk}{x^j}'{x^k}'+C^i_j{x^j}'
+D^i=0, \quad i,j,k,l=2,\ldots,n, \label{(r5)}
\end{equation}
provided 
\begin{eqnarray}
&&\Delta^i_{klm}=J^i_jA^j_{klm}+G^i_{mj}B^j_{kl},\nonumber\\
&&\Lambda^i_{kl}=J^i_jB^j_{kl}+G^i_{lj}C^j_k,\nonumber\\
&&\Omega^i_k=J^i_jC^j_k+G^i_{kj}D^j,\nonumber\\
&&E^i=J^i_jD^j,\nonumber\\
&&G^i_{pj}A^j_{klm}=0.\label{(r55)}
\end{eqnarray}

The relations (\ref{(r55)}) can be obtained by solving for the second derivative in terms of the 
first order derivatives and inserting these into equation (\ref{(r1)}). The last equation of (\ref{(r55)})
tells us that not all the $A^i_{jklm}$ coefficients are independent. As a matter of fact if we replace 
these by $A_{kl}$ in (\ref{(r5)}), then it turns out that this relation
in (\ref{(r55)}) will now be identically satisfied. What transpires is that the quartic term disappears automatically
due to $G^i_{pj}$ being skew symmetric in the lower indices and ${x^p}'{x^j}'$ appearing symmetrically. 
One also needs then to adjust
the relation (\ref{(r55)}a) in the latter case by
\begin{equation}
\Delta^i_{klm}=J^i_kA_{lm}+G^i_{mj}B^j_{kl}.\label{(r57)}
\end{equation}
The remaining equations of (\ref{(r55)}) are the same.

There are two branches of the linearization problem by point transformations for a system of $n-1$ second-order
ODEs. One is the general form (\ref{(r1)}) owing to the arbirariness of the $\Delta^i_{jkl}$ coefficients. 
The other is the form (\ref{(14)}) in which the cubic coefficients are fewer in number. In the case of two second-order ODEs,
equations (\ref{(50)}), we have obtained explicit linearization criteria as encapsulated in Theorem 2
and their corollaries.

In the general equation (\ref{(r1)}) there are 
$(n-1)n(n^2+6n-1)/6$ coefficients while for (\ref{(14)}) there are $(n-1)n(n+2)/2$ independent coefficients.
It would be of interest to find practical criteria for the reduction of equation (\ref{(r1)}) to the simplest system via point transformations for $n=3$. Of course it is of great interest to do this for the general system (\ref{(r1)}) for
$n\ge4$.

If one has a system of the form (\ref{(r1)}) with known coefficients which is reducible to the free particle system (\ref{(r4)}) by point transformations, then one can utilise (\ref{(r2)}) to construct a linearizing point transformation. Also, we can obtain linearizing 
point transformations for system (\ref{(14)}), if it is linearizable to the simplest system (\ref{(r4)}), by invoking (\ref{(r2)})
together with (\ref{(r55)}).

In particular, one can find linearizing point transformations for the system (\ref{(50)}) in a similar manner
by solving the system (\ref{(r55)}).

Instead of using the system (\ref{(r2)}) in order to construct a linearizing point transformation 
there are other ways as pointed out earlier. One is to go to the higher space, once one has the coefficients at hand, and use 
(\ref{(12)}) for which $g_{ij}({\bf u})$ must be the identity matrix and where we may set $u^1$ to be the independent
variable. Yet a third approach is that of mapping symmetry generators of the linearizable system, if known, to the free particle
generators.

\section{Examples}
We present examples to illustrate our results. We have 
$y$ and $z$ as the dependent variables. Also the $'$ below denotes differentiation
with respect to $x$. Moreover, we have included one example that does not satisfy our linearization criteria but belongs to the more general class (\ref{(r1)}) which is linearizable.

{\bf 1.} Consider the anisotropic oscillator system
\begin{eqnarray}
y''+\omega_1(x)y=0,\nonumber\\
z''+\omega_2(x)z=0,\label{(56)}
\end{eqnarray}
The coefficients of system (\ref{(56)}) satisfy the conditions
(\ref{(55)}) provided $\omega_1=\omega_2$. Hence in order for the
system (\ref{(56)}) to be reducible to the free particle system one
must have isotropy.

{\bf 2.} The simple linear system
\begin{eqnarray}
y''+z=0,\nonumber\\
z''+z=0,\label{(57)}
\end{eqnarray}
do not satisfy the conditions (\ref{(55)}). Thus this system is not
transformable pointwise to the free particle system. This system
does not have a Lagrangian formulation as well \cite{dou}.

{\bf 3.} For the quadratic system
\begin{eqnarray}
y''-y'+y'^2=0,\nonumber\\
z''-z'+z'^2=0,\label{(58)}
\end{eqnarray}
all conditions (\ref{(53)}) are satisfied. Therefore the system
(\ref{(58)}) is reducible to the simplest system. A point
transformation that does the job is
\begin{equation}
u=\exp x,\; v=\exp y,\; w=\exp z,\label{(l1)}
\end{equation}
where $u$ is the independent variable. This can be constructed by
going to the higher space as we have illustrated for the scalar ODEs
in section 3.

{\bf 4.} Consider the cubically semi-linear system
\begin{eqnarray}
y''+\frac1{x}y'+y'^2+(\frac{x}{y}+{x\over y^2})y'^3=0,\nonumber\\
z''+\frac1{x}z'+z'^2+2y'z'+(\frac{x}{y}+{x\over
y^2})y'^2z'=0,\label{(59)}
\end{eqnarray}
For the system (\ref{(59)}) all the conditions (\ref{(51)}) hold. A
linearizing point transformation to the simplest system is
\begin{equation}
u=\ln xy,\; v=\exp y,\; w=\exp(y+z),\label{(l2)}
\end{equation}
in which $u$ is the independent variable.

{\bf 5.} Finally the system
\begin{eqnarray}
4yz^2y'^2+4y^2zy'z'+2xz^2y'^3+8xyzy'^2z'+2xy^2y'z'^2+2xy^2zy'z''=y^2z^2y''+2xy^2zz'y'',\nonumber\\
y''+xzy'y''+xyz'y''-xz^2y'^2y''-xyzy'z'y''=y'(zy'+yz')(zy'+yz'+2xy'z'+xyz''),\label{(60)}
\end{eqnarray}
is not of the form given in Theorem 2. It is of the form given in (\ref{(r1)}) and is linearizable 
by means of the point transformation
\begin{equation}
u=x\exp(yz),\; v=xy^2z^2,\; w=y,\label{(l3)}
\end{equation}
where $u$ is the independent variable.

\section{Concluding remarks}

Aminova and Aminov \cite{ami} had provided a procedure of projecting
down 1 dimension from a system of $n$ geodesic equations to $n-1$
cubically semi-linear ODEs. Separately, we had provided \cite{mah4}
linearizability criteria for a square quadratically semi-linear
system. These were used together to derive linearizability criteria
for a single cubically semi-linear equation by projecting down from
a system of 2 quadratically semi-linear equations. This provided an
alternate method to prove Lie's general result for linearizability
of a single non-linear equation. It led naturally to an extension of
the linearization criteria via point transformations from a scalar
second-order ODE as obtained by Lie \cite{lie} to a system of two
cubically semi-linear ODEs of the form (\ref{(50)}). These provided
necessary and sufficient conditions for reduction to the simplest
system and hence $sl(4,\R)$ symmetry algebra for equations of the
form (\ref{(50)}). Moreover, Theorem 2 provides criteria for the
reduction of linear systems of two equations to the free particle
system.

Lie had demonstrated \cite{lie} that only cubically semi-linear
scalar equations of order two are linearizable in general. 
As such, it could have been
hoped that the projection procedure will provide the complete
solution of the linearizability problem for the system of 2
non-linear ODEs. That hope is doomed from the start as there are 5
classes of systems of 2 cubically non-linear equations that are linearizable
by point transformations, having
different symmetry algebras. Moreover, the maximum symmetry algebra class of such systems
of two equations is one branch of the linearization problem
via point transformations as the general class is represented by (\ref{(r1)}).
Why do we get a unique class in the
former case and 5 in the latter? Furthermore, how many distinct
classes should there be for a system of $n$ cubically semi-linear
ODEs?

We start by noting that the projection procedure and linearizability can be equally well
adopted for an arbitrary system of $n$ quadratically semi-linear
second order ODEs reduced to $n-1$ cubically semi-linear second
order ODEs. There are two branches for the linearization problem for systems admitting
the maximal algebra for $n\ge3$. There is enormous computational
complications that arise. As such, one would need an algebraic
computational code to deal with larger systems. A code has, indeed, been prepared to construct
the metric coefficients given the Christoffel symbols \cite{fred}. That can 
be extended to deal with the linearization of larger systems. Now observe that in
projecting down from the system of $n$ dependent variables to $n-1$ variables,
the Christoffel symbols are reduced from $n^2(n+1)/2$ by $n$, to
give $(n-1)n(n+2)/2$ independent coefficients. Since we now have
$n-1$ equations, each with its own cubic function, there are
$(n-1)n/2$ cubic coefficients for the reduced system. If the number of
coefficients left over after losing $n$ equals the number of
coefficients of the reduced system, we can determine one set of
coefficients in terms of the other. The two expressions are
obviously equal for $n=2$ and the former is greater than the latter
for $n>2$. As such, the coefficients of the cubic system can be
determined uniquely in terms of the quadratic system for $n=2$, i.e.
for a scalar cubically semi-linear system. For larger systems there
will be infinitely many ways to write the former in terms of the
latter. Hence there is a unique solution to the linearizability
problem only for the scalar cubically semi-linear equation and many
solutions for systems of cubically semi-linear systems!

The second question remains and has, in fact, been compounded. It is
known that there are 5 and not infinitely many distinct classes.
Why? The point is that all distinct ways of writing the cubic system
coefficients in terms of the quadratic system coefficients will not
give independent criteria as there will be transformations
permissible from one definition to another. The point is to
determine those that are distinct. Another way of looking at what we
have done is to note that we have asked that the original system
correspond to a system of geodesic equations in flat space. Then the
projection gives the reduced system, which must also be of geodesics
in an $(n-1)$-dimensional flat space. Even if the original geodesics
were curved, the projected geodesics could correspond to straight
lines. For example, if the original space was a sphere and one
projects along the plane containing the geodesic to a plane
perpendicular to it, the resulting projected curve would be a
straight line.

The minimal dimension for a system of $n$ second-order ODEs to be linearizable by point transformation
is $2n+1$. The maximum dimension of the symmetry algebra is $(n+1)(n+3)$ which corresponds to $sl(n+2,\R)$. 
The other submaximal 
symmetry alebras besides that of dimension $2n+1$ range from $2n+2$ to $(n+2)^2/2$ for $n$ even
and $[(n+2)^2+1]/2$ for $n$ odd. Thus for $n=2$ we have the mimimum dimension to be 5 and other submaximal
algebra dimensions are 6, 7 and 8. The maximum dimension for $n=2$ is 15. For $n=3$ the minimum dimension is 7
and the next to maximum is 13. The maximum is 24. Thus for this case there are 8 classes. Generally, for $n=2m$, 
the number of classes is $2m^2+3$ and for $n=2m-1$ it is $2m^2-2m+4$.

It would be important to find ways of providing the linearizability
criteria for the cases of the other symmetry algebras.

\section*{Appendix}

We take $j=1$ in the third set of (\ref{(25)}) as the 3 dependent
equations and discard them. The invocation of the first set of 9
equations of (\ref{(25)}) gives
\begin{eqnarray}
&&\frac12C^3_{2x}-D^3_y+\frac14C^3_3C^3_2+\frac14C^2_2C^3_2-D^2B^3_{22}-D^3B^3_{23}=0,\nonumber\\
&&B^3_{22x}-\frac12C^3_{2y}-A_{22}D^3+\frac12C^3_2B^2_{22}+\frac12C^3_3B^3_{22}-\frac12C^2_2B^3_{22}-\frac12B^3_{23}C^3_2=0,\nonumber\\
&&\Gamma_{12,y}^1=-A_{22x}-A_{22}\Gamma^2_{12}+C^2_2A_{22}+\Gamma^1_{12}B^2_{22}+\Gamma^1_{12}\Gamma_{12}^1+\frac12B^3_{22}\Gamma^3_{33}
-\frac12B^3_{22}B^3_{33}+\frac12 C^3_2A_{23},\nonumber\\
&&\Gamma_{12,x}^2=D_y^2+D^2\Gamma^1_{12}+\Gamma^2_{12}\Gamma^2_{12}-\frac14C^2_3C_2^3-\Gamma^2_{12}C^2_2+B^2_{22}D^2
+D^3B^2_{23}-\frac12D^3B^3_{33}+\frac12D^3\Gamma^3_{33},\nonumber\\
&&\Gamma^2_{12,y}=-\frac13B^2_{22x}+\frac23C^2_{2y}+\Gamma^1_{12}\Gamma^2_{12}+\frac14C^3_2\Gamma^3_{33}-\frac14C^3_2B^3_{33}\nonumber\\
&&+D^2A_{22}+\frac23D^3A_{23}-\frac16C^2_3B^3_{22}+\frac1{6}B^2_{23}C^3_2,\nonumber\\
&&\Gamma^1_{12,x}=-\frac23B^2_{22x}+\frac13C^2_{2y}+\Gamma^1_{12}\Gamma^2_{12}+\frac14C^3_2\Gamma^3_{33}-\frac14C^3_2B^3_{33}
+D^2A_{22}\nonumber\\
&&+\frac13D^3A_{23}-\frac13C^2_3B^3_{22}+\frac13B^2_{23}C^3_2\nonumber\\
&&B^3_{23x}-\frac13B^2_{22x}+\frac16C^2_{2y}-\frac43D^3A_{23}-\frac23B^3_{22}C^2_3+\frac2{3}B^2_{23}C^3_2-\frac12C^3_{3y}=0,\nonumber\\
&&\Gamma_{33,y}^3=-2A_{23x}+B^3_{33y}-2\Gamma^2_{12}A_{23}+A_{23}C^2_2+2\Gamma_{12}^1B^2_{23}+\Gamma_{12}^1\Gamma_{33}^3-\Gamma_{12}^1B^3_{33}
\nonumber\\
&&+\Gamma_{33}^3B^3_{23}-B^3_{33}B^3_{23}+A_{22}C^2_3+A_{23}C^3_3\nonumber\\
&&\Gamma_{33,x}^3=-2B^2_{23x}+B^3_{33x}+C^2_{3y}+2D^2A_{23}-C^2_3B^3_{23}+C^2_3\Gamma_{12}^1-\Gamma^2_{12}B^3_{33}+C^2_3B^2_{22}+B^2_{23}C^3_3\nonumber\\
&&-B^2_{23}C^2_2-\frac12C^3_3B^3_{33}+\frac12B^3_{33}C^2_2+\frac12C^3_3\Gamma^3_{33}-\frac12C^2_2\Gamma^3_{33}+\Gamma^3_{33}\Gamma^2_{12}.
\label{(A1)}
\end{eqnarray}
The second set of 9 equations of (\ref{(25)}) yields
\begin{eqnarray}
&&\Gamma_{12,x}^1=-B^3_{23x}+\frac12C^3_{2z}+D^3A_{23}-\frac12C^3_2B^2_{23}+\frac14C^3_2B^3_{33}+\frac14C^3_2\Gamma_{33}^3\nonumber\\
&&-\frac12C^3_3B^3_{23}+\frac12C^2_2B^3_{23}+\Gamma_{12}^2\Gamma_{12}^1,\nonumber\\
&&\Gamma_{12,z}^1=-A_{23x}-A_{23}\Gamma_{12}^2+A_{23}C^2_2+\Gamma_{12}^1B^2_{23}+\frac12\Gamma_{33}^3B^3_{23}\nonumber\\
&&\frac12\Gamma_{33}^3\Gamma_{12}^1-\frac12B^3_{33}B^3_{23}-\frac12B^3_{33}\Gamma_{12}^1+\frac12A_{33}C^3_2,\nonumber\\
&&\Gamma_{12,x}^2=-\frac12C^3_{3x}+\frac12C^2_{2x}+D^3_z+\frac12\Gamma_{33}^3D^3+\frac12D^3B^3_{33}-\frac14C^3_2C^2_3-\frac14C^3_3C^3_3\nonumber\\
&&+\frac14C^2_2C^2_2+\Gamma_{12}^2\Gamma_{12}^2-C^2_2\Gamma_{12}^2+B^2_{23}D^2+\Gamma_{12}^1D^2,\nonumber\\
&&\Gamma_{12,z}^2=-B^2_{23x}+2A_{23}D^2-\frac12C^2_3B^3_{23}+\frac12B^2_{33}C^3_2+C^2_{2z}+\frac12C^2_3\Gamma_{12}^1+\frac14C^3_3\Gamma_{33}^3
\nonumber\\
&&-\frac14C^2_2\Gamma_{33}^3+\frac12\Gamma_{33}^3\Gamma_{12}^2-\frac14B^3_{33}C^3_3+\frac14C^2_2B^3_{33}-\frac12B^3_{33}\Gamma^2_{12}+A_{33}D^3,
\nonumber\\
&&\frac12C^2_{3x}-D^2_z+\frac14C^2_3C^3_3+\frac14C^2_3C^2_2-B^2_{23}D^2-B^2_{33}D^3=0,\nonumber\\
&&B^2_{33x}-\frac12C^2_{3z}-D^2A_{33}+\frac12C^2_3B^3_{33}-\frac12B^2_{23}C^2_3-\frac12B^2_{33}C^3_3+\frac12B^2_{33}C^2_2=0,\nonumber\\
&&\Gamma_{33,x}^3=B^3_{33x}-4B^2_{23x}+6A_{23}D^2-2C^2_3B^3_{23}+2B^2_{33}C^3_2+2C^2_{2z}+C^2_3\Gamma_{12}^1+\frac12C^3_3\Gamma_{33}^3\nonumber\\
&&-\frac12C^2_2\Gamma_{33}^3+\Gamma_{33}^3\Gamma_{12}^2-\frac12B^3_{33}C^3_3+\frac12C^2_2B^3_{33}-B^3_{33}\Gamma_{12}^2+2A_{33}D^3,\nonumber\\
&&\Gamma_{33,x}^3=\frac12C^3_{3z}+\frac12C^2_{2z}-B^2_{23x}+2A_{23}D^2+C^2_3\Gamma_{12}^1+\frac12C^3_3\Gamma_{33}^3-\frac12C^2_2\Gamma_{33}^3\nonumber\\
&&+\Gamma_{12}^2\Gamma_{33}^3-\frac12C^3_3B^3_{33}+\frac12C^2_2B^3_{33}-B^3_{33}\Gamma_{12}^2+2A_{33}D^3,\nonumber\\
&&\Gamma_{33,z}^3=-2A_{33x}+B^3_{33z}-2A_{33}\Gamma_{12}^2+C^2_2A_{33}+2\Gamma_{12}^1B^2_{33}+\frac12\Gamma^3_{33}\Gamma_{33}^3
-\frac12B^3_{33}B^3_{33}\nonumber\\
&&+A_{23}C^2_3+A_{33}C^3_3,\label{(A2)}
\end{eqnarray}
The last set of the 9 equations of (\ref{(25)}) result in 6 independent conditions
\begin{eqnarray}
&&\Gamma_{12,y}^1=-B^3_{23y}+B^3_{22z}+\frac12C^3_2A_{23}-B^3_{22}B^2_{23}+\frac12B^3_{22}B^3_{33}+\frac12B^3_{22}\Gamma_{33}^3\nonumber\\
&&-B^3_{23}B^3_{23}+\Gamma_{12}^1\Gamma_{12}^1-\frac12C^3_3A_{22}+\frac12C^2_2A_{22}-A_{22}\Gamma_{12}^2+B^3_{23}B^2_{22}+
B^2_{22}\Gamma_{12}^1,\nonumber\\
&&\Gamma_{12,z}^1=\frac13B^3_{23z}+\frac16C^3_2A_{33}-B^2_{33}B^3_{22}-\frac1{6}C^3_3A_{23}+\frac1{6}C^2_2A_{23}-A_{23}\Gamma_{12}^2
+B^2_{23}B^3_{23}-\frac12B^3_{23}B^3_{33}\nonumber\\
&&+\frac12B^3_{23}\Gamma_{33}^3+B^2_{23}\Gamma_{12}-\frac12B^3_{33}\Gamma^1_{12}+\frac12\Gamma_{12}^1\Gamma_{33}^3
+\frac23B^2_{23y}-\frac13B^3_{33y}-\frac23B^2_{22z}+\frac13C^2_3A_{22},\nonumber\\
&&\Gamma_{33,y}^3=\frac43B^3_{23z}+\frac23C^3_2A_{33}-2B^2_{33}B^3_{22}-\frac23C^3_3A_{23}+\frac23C^2_2A_{23}
-2A_{23}\Gamma_{12}^2 +2B^3_{23}B^2_{23}\nonumber\\
&&-B^3_{23}B^3_{33}+B^3_{23}\Gamma_{33}^3+2B^2_{23}\Gamma_{12}^1-B^3_{33}\Gamma_{12}^1
+\Gamma_{12}^1\Gamma_{33}^3+\frac23B^2_{23y}-\frac13B^3_{33y}-\frac23B^2_{22z}+\frac13C^2_3A_{22},\nonumber\\
&&\Gamma_{33,z}^3=2B^2_{33y}-2B^2_{23z}+B^3_{33z}-2A_{33}\Gamma_{12}^2+2B^2_{22}B^2_{33}+2B^2_{33}\Gamma_{12}^1+\frac12\Gamma_{33}^3\Gamma_{33}^3\nonumber\\
&&+C^2_3A_{23}-2B^2_{23}B^2_{23}+2B^2_{23}B^3_{33}-\frac12B^3_{33}B^3_{33}-2B^2_{33}B^3_{23},\nonumber\\
&&-A_{23y}+A_{22z}-A_{22}B^2_{23}-A_{23}B^3_{23}+A_{23}B^2_{22}+A_{33}B^3_{22}=0,\nonumber\\
&&-A_{33y}+A_{23z}-A_{22}B^2_{33}-A_{23}B^3_{33}+A_{23}B^2_{23}+A_{33}B^3_{23}=0,\label{(A3)}
\end{eqnarray}

\section*{Acknowledgements}
FM thanks the HEC of Pakistan for granting him a visiting professor
and NUST-CAMP for hospitality during which time this work was
undertaken. AQ is grateful to the DECMA Centre and the School of Computational
and Applied Mathematics for hosting him when this work was completed.

\end{document}